\documentclass[10pt]{article}

\usepackage{amsmath, amsthm, amsfonts}
\usepackage{graphicx}
\usepackage{verbatim}

\newtheorem{theorem}{Theorem}[section]

\newtheorem{question}[theorem]{Question}

\begin{document}

\begin{center}
\Large \bf Convergence properties of Donaldson's $T$-iterations on the Riemann sphere
\end{center}

\begin{center}
Morgan Sherman\footnote{%
    Part of this work was carried out while the author was visiting the Mathematics Department of Harvard University, and he wishes to thank them for their hospitality.
}\\
California Polytechnic State University
\end{center}

\smallskip

\begin{center}
\begin{minipage}{\textwidth}
\small
{\bf Abstract.}  In \cite{donaldson5} Donaldson gives three operators on a space of Hermitian metrics on a complex projective manifold: $T, T_{\nu}, T_K.$  Iterations of these operators converge to {\it balanced} metrics, and these themselves approximate constant scalar curvature metrics.  In this paper we investigate the convergence properties of these iterations by examining the case of the Riemann sphere as well as higher dimensional $\mathbb{CP}^n$.
\end{minipage}
\end{center}

\section{Introduction}

Let $X$ be a compact complex manifold with a positive holomorphic line bundle $L$.
A long-standing open problem in K\"ahler geometry, building on Yau's solution of the Calabi conjecture \cite{yau1}, is to find sufficient conditions for the existence of a constant scalar curvature K\"ahler metric  in $c_1(L)$.  Another is as follows: can such a metric be obtained naturally as a limit of algebraic metrics via embeddings of $X$ into $\mathbb{P} \mathrm{H}^0 (X, L^k)$?

This idea of approximating K\"ahler metrics by restricting Fubini-Study metrics, advocated by Yau over the years, has led to the development of a rich theory relating analysis and notions of stability in the sense of geometry invariant theory (see \cite{yau2}, \cite{tian1}, \cite{tian2}, \cite{donaldson3}).  In a fundamental paper, Donaldson \cite{donaldson2} showed that, under an assumption on the space of automorphisms, the metrics induced from \emph{balanced} embeddings \cite{zhang}  of $X$ into projective space by sections of $L^k$ converge to the constant scalar curvature metric as $k\rightarrow \infty$.  The balanced condition means that
\[
        \int_X \frac{Z_i \overline{Z_j}}{|Z|^2}\, d\mu \ = \ c\delta_{ij},
\]
(where $d\mu$ is the volume form on $X$ induced by the Fubini-Study metric and $c$ is a constant depending on the data $(X, L^k)$ and not on the particular embedding) and this is equivalent to the Chow stability of the embedding \cite{zhang}, \cite{luo}, \cite{ps1}.

Recently \cite{donaldson4, donaldson5}, Donaldson has devised iterative procedures on the space of Hermitian metrics on $\mathrm{H}^0(X, L^k)$  to find approximations to these balanced metrics.  For sufficiently large $k$, these approximations are close to a constant scalar curvature metric.  Explicit numerical computations, focused in the Calabi-Yau case where there are possible applications to string theory, have been obtained in \cite{donaldson5}, \cite{dklr1}, \cite{dklr2}.  See also \cite{hw}, \cite{keller2}, and \cite{rubinstein} where different methods are used.

Donaldson's three iterative maps $T$, $T_{\nu}$, and $T_K$, described below, are interesting in their own right. Indeed, as pointed out in \cite{donaldson5}, it is likely that these maps can be viewed as discrete approximations to the Ricci and Calabi flows.

Instead of pursuing general questions of existence, in this paper we pick a simple compact complex manifold -- the Riemann sphere -- and investigate the convergence properties of each of $T$, $T_{\nu}$, and $T_K$ on the space of Hermitian metrics induced from Veronese embedings into $\mathbb{CP}^n.$  In section \ref{section: Higher dimensional projective space} we briefly investigate the case $\mathbb{CP}^n$ when $n>1$.

There is a natural notion of distance  on the space of Hermitian metrics $\mathrm{GL}(n+1, \mathbb{C})/\mathrm{U}(n+1)$, and indeed as $k$ increases  this distance function is expected \cite{ps2}  to approximate that on the infinite dimensional space of K\"ahler metrics  \cite{mabuchi1}, \cite{semmes}, \cite{donaldson1}, \cite{chen}.  A natural question one might ask is:  Are the $T, T_{\nu},$ or  $T_K$ iterations distance reducing on the space of metrics?  In section \ref{subsection: The effect on distance} we show that the $T$ operator does \emph{not} satisfy this property.

One goal of this study was to find an effective bound on the distance between the $r$th iteration of a metric under $T$, $T_{\nu}$, or $T_K$ and the limiting balanced metric.  One is proposed in section
\ref{subsection: The effect on distance}.  In section \ref{subsection: Asymptotic behavior} we list the observed asymptotic behavior of each of these iterations.
In section \ref{section: Examples} we give some examples.  In section
\ref{section: Higher dimensional projective space} we investigate the case for higher dimensional projective space.

It has recently come to the author's attention that on Julien Keller's web site \cite{keller} one can find a program to compute a Ricci flat metric on a particular $K3$ surface using the techniques of Donaldson on which this paper is based.  More information can be found there.  All computations and all graphs in this paper were done using the software Maple 9.
\bigskip

\noindent
{\bf Acknowledgements.}  The author is grateful to Ben Weinkove for introducing him to this problem, and for answering endless questions.  This paper would not have been possible without his help.  The author would also like to thank the referee for many useful comments and suggestions which helped to improve this paper.


\section{The $T$, $T_{\nu}$, and $T_K$ operators}\label{section: The T operators}

Let $X$ be an $n$ dimensional complex projective manifold, and $L \rightarrow X$ an ample line bundle.  In \cite{donaldson5} Donaldson examines three different actions on the space of Hermitian metrics on $\mathrm{H}^0(X, L^k)$:  $T, T_{\nu}, T_K.$
We briefly recall how he defines each.

Given a Hermitian metric $G$ on $\mathrm{H}^0(X, L^k)$ and an orthonormal basis $\{s_i\}$ with respect to $G$, one defines the Fubini-Study metric $h= \mathrm{FS}(G)$ on the line bundle $L^k$ by the requirement that $\sum_i |s_i|_{h}^2 = 1$.  The result is independent of the orthonormal basis chosen.  Now given this metric $h$ on $L^k$ we define a new Hermitian metric on $\mathrm{H}^0(X, L^k)$, denoted $\mathrm{Hilb}(h)$, by
\[
    \| s \|_{\mathrm{Hilb}}^2 = R\int_{X} |s|_h^2 \omega_h^n/n!
\]
where $\omega_h$ is the K\"ahler form $-\sqrt{-1}\partial\overline{\partial}\log h$ and
where $R$ is the constant
\[
    R = \frac{\mathrm{dim}\,  \mathrm{H}^0(X, L^k)}{\mathrm{Vol}(X, \omega_h^{n}/n!)}.
\]
This defines the $T$ map:  $T(G)=\mathrm{Hilb}(\mathrm{FS}(G))$.

The $T_{\nu}$ map is defined analogously, but instead of the volume form $\omega_h^n/n!$ we fix a volume form $\nu$ of our choosing.  As above we set
\[
    \| s \|_{\mathrm{Hilb}_{\nu}} = R_{\nu} \int_X |s|_h^2 \nu,
\]
where
\[
    R_{\nu} = \frac{\mathrm{dim}\,  \mathrm{H}^0(X, L^k)}{\mathrm{Vol}(X, \nu)}.
\]
Then we define $T_{\nu}(G) = \mathrm{Hilb}_{\nu}(\mathrm{FS}(G)).$

The $T_K$ function is defined in case $L^k = K^{-p}$, where $K$ is the canonical bundle.  Again we only modify the volume form, this time choosing
\[
    \omega_{G,K} = \left( \sum s_i \otimes \overline{s}_i \right)^{-1/p}.
\]
The resulting metric on $\mathrm{H}^0(X, L^k) = \mathrm{H}^0(X, K^{-p})$ is given as above:
\[
    \| s \|_{\mathrm{Hilb}_K} = R_K \int_X |s|_h^2 \omega_{G,K}
\]
where
\[
    R_K = \frac{\mathrm{dim}\,  \mathrm{H}^0(X, L^k)}{\mathrm{Vol}(X, \omega_{G,K})}.
\]
As before set $T_K(G) = \mathrm{Hilb}_K(\mathrm{FS}(G)).$

A Hermitian metric $G$ is {\it balanced with respect to } $T$ (resp. $T_{\nu}, T_K$) if $T(G) = G$ (resp. $T_{\nu}(G) = G, T_K(G)=G$).  The basic philosophy is that if $F = T, T_{\nu}, T_K$ and if there exists some balanced metric, then starting with any Hermitian metric $G$ the iterations $F^{(r)}(G)$ should tend to a balanced metric as $r$ tends to infinity (see \cite{donaldson5} and also \cite{sano}).   In this paper we will concern ourselves only with a very simple case and study in some detail the properties of this convergence.

Specifically we take as our manifold the Riemann sphere $X = \mathbb{CP}^1$ and line bundle $L = O_X(1)$.  We note that the presence of the automorphism group $\mathrm{SL}(2, \mathbb{C})$ means that, strictly speaking, some aspects of the theory may need to be developed further, in the manner of \cite{mabuchi2} for example, but since we are focusing on numerical results here we will not dwell on this issue.  Fix a holomorphic coordinate $z \in \mathbb{C}.$  Then $\mathrm{H}^0(X,L^k) = \mathrm{H}^0(\mathbb{CP}^1, O(k)) \cong \mathbb{C}^{k+1}$ has basis $1, z, z^2, \ldots, z^k$.  Hermitian metrics can now be associated with $(k+1)\times(k+1)$ positive definite Hermitian matrices.  For the $T_{\nu}$ function we fix our volume form $\nu$ as the standard Fubini-Study form
\begin{equation} \label{nu form}
    \nu = \sqrt{-1}\partial\overline{\partial}\log(1+|z|^2)
    = \frac{\sqrt{-1}}{(1+|z|^2)^2} \mathrm{d}z\wedge\mathrm{d}\overline{z}.
\end{equation}
In the case of the $T_K$ map we note that $K = O(-2)$, hence $L^k = K^{-p}$ precisely when $k=2p.$

We simplify further by considering only those metrics invariant under the $S^1$ action $z \mapsto e^{i\theta}z$ on the Riemann sphere.  This restricts our attention to {\it diagonal} positive definite Hermitian $(k+1)\times(k+1)$ matrices $G$.  We will suppose $G$ has entries $a_0^{-1}, a_1^{-1}, \ldots, a_k^{-1}$ -- taking inverses simplifies later computations -- and we will use the notation
\[
    G = \left( a_0, a_1, \ldots, a_k \right)
\]
to denote this metric.  Each of $T, T_{\nu},$ and $T_K$ is a function of $(a_0, a_1, \ldots, a_k)$, and in the remainder of this section we write them down explicitly.

We begin with $T.$  Taking $G$ as above we can pick the orthonormal basis $\{s_i = \sqrt{a_i}\,z^i,\ i=0,\ldots,k\}.$  Then
\[
    h = FS(G) = \left( \sum a_i |z|^{2i} \right)^{-1}
\]
and we calculate
\[
    \omega_h = \sqrt{-1}\partial\overline{\partial}\log\left(\sum a_i|z|^{2i} \right)
        = \sqrt{-1} \frac{\sum_{i>j}a_ia_j(i-j)^2|z|^{2(i+j-1)}}{\left(\sum a_i|z|^{2i}\right)^2}
            \mathrm{d}z\wedge\mathrm{d}\overline{z}.
\]
Write $T(a_0, \ldots, a_k) = (\tilde{a}_0, \ldots, \tilde{a}_k).$  Then
\[
    \tilde{a}_q^{-1} = R \int_{\mathbb{C}} |z|^{2q} h \omega_h,
\]
where $R = (k+1)/\mathrm{Vol}(X, \omega_h)$.
Using polar coordinates $z=re^{i\theta}$ and setting $x = r^2$ we get
\[
    \tilde{a}_q = 1 / \left(
    2\pi R \int_0^{\infty} \frac{ \sum_{i>j}a_ia_j(i-j)^2x^{i+j-1} }{ \left( \sum a_i x^{i} \right)^2}
        x^q\mathrm{d}x
    \right).
\]
Thus after substituting for $R$ we find
\begin{equation} \label{T map}
    T : a_q \mapsto
    \frac{
    \int_0^{\infty} \frac{ \sum_{i>j}a_ia_j(i-j)^2x^{i+j-1} }{ \left( \sum a_i x^{i} \right)^2}
        \mathrm{d}x
        }{
    (k+1) \int_0^{\infty} \frac{ \sum_{i>j}a_ia_j(i-j)^2x^{i+j-1} }{ \left( \sum a_i x^{i} \right)^3}
        x^q\mathrm{d}x}, \ q=0,1,\ldots,k.
\end{equation}

By a similar computation, noting the $T_{\nu}$ map has the simpler volume form (\ref{nu form}), we find
\begin{equation} \label{Tnu map}
    T_{\nu}: a_q \mapsto
    \left( (k+1) \int_0^{\infty} \frac{x^q \mathrm{d}x}{(1+x)^2 \sum a_i x^i} \right)^{-1}
    ,\ q=0,1,\ldots,k.
\end{equation}

For the $T_K$ map the volume form is
\[
    \omega_{G,K} = \sqrt{-1}\left( \sum a_i |z|^{2i} \right)^{-1/p} \mathrm{d}z\wedge\mathrm{d}\overline{z}
\]
and we calculate as above:
\begin{equation} \label{TK map}
    T_K : a_q \mapsto
    \frac{
    \int_0^{\infty} \left( \sum a_i x^{i} \right)^{-2/k} \mathrm{d}x
    }{
    (k+1) \int_0^{\infty} \left( \sum a_i x^{i} \right)^{-1-2/k} x^q\mathrm{d}x
    },\ q=0,1,\ldots,k.
\end{equation}

Often it is simpler still to work with ($S^1$-invariant) metrics invariant under the inversion $z \mapsto z^{-1}.$  We call such metrics \emph{palindromic} as they are characterized as those metrics $(a_0, a_1, \ldots, a_k)$ which satisfy
\[
    a_0 = a_k, \ a_1 = a_{k-1}, \ldots,\ a_{\lfloor k/2 \rfloor} = a_{\lceil k/2 \rceil}.
\]
Thus in the palindromic case there are exactly $\lceil k/2 \rceil$ real (positive) parameters, while in the non-palindromic case there are $k+1.$  However we note that for any of the operators $F = T, T_{\nu}, T_K$, and any starting metric $(a_0, \ldots, a_k)$, if we let $(\tilde{a}_0, \ldots, \tilde{a}_k)$ denote the metric after an application of $F$, then we have a relation
\begin{equation}\label{eqn:  metric relation}
    \sum_{i=0}^k \frac{a_i}{\tilde{a}_i} = k+1.
\end{equation}
This is immediately verified by checking formulas (\ref{T map}), (\ref{Tnu map}), (\ref{TK map}).


\section{Findings}\label{section: Findings}

In investigating the behavior of the convergence of a sequence of Hermitian metrics we need to decide what we mean when we say two metrics are close.
Let $M = \mathrm{GL}(k+1, \mathbb{C})/U(k+1)$ be the space of Hermitian metrics on $\mathbb{CP}^k$.  The $\mathrm{GL}(k+1,\mathbb{C})$-invariant K\"ahler metric is given by the form $g_H(U,V) = tr(H^{-2}UV)$ where $U, V$ are in the tangent space to $H$ on $M$.  Geodesics on $M$ are given by the images of one-parameter subgroups, e.g.
\[
    \left( \begin{array}{ccc} e^{\alpha_0 t} \\ & \ddots \\ & & e^{\alpha_k t} \end{array} \right).
\]
Let $A = (a_0, \dots, a_k)$, and $B = (b_0, \ldots, b_k)$ be two metrics in $M$.  Writing $a_i = e^{\alpha_i}$ and $b_i = e^{\beta_i}$ for $i=0, \ldots, k$ we find the geodesic from $A$ to $B$ is given by $P(t), 0 \leq t \leq 1$, where $P(t)$ is the diagonal matrix with entries $e^{(\beta_i-\alpha_i) t + \alpha_i}, i=0,\ldots,k$.
Now we can calculate the distance between $A$ and $B$ as
$\int_0^1 \left| \frac{dP}{dt} \right|_P \, dt = \sqrt{ \sum (\beta_i - \alpha_i)^2 }$ or
\begin{equation}\label{eqn:  distance formula}
    \mathrm{dist}(A,B)  =
    \sqrt{ \sum_{i=0}^k \left( \log \frac{b_i}{a_i} \right)^2 }\ .
\end{equation}

One goal is then to understand how well the $r$th iteration of $F = T, T_{\nu}, T_K$ applied to a Hermitian metric $G$ approximates the limiting balanced metric $B := F^{(\infty)}(G).$  That is we wish to understand the function
\[
    \mathrm{err}_{F, k} (G, r) = \mathrm{dist}\left( F^{(r)}(G),\ F^{(\infty)}(G) \right).
\]
In particular we would like to give an effective bound:
\[
    \mathrm{err}_{F, k} (G, r) < \mathrm{bnd}_{F, k} (d, r)
\]
where $d = \mathrm{dist}(G,B).$  We propose such a bound in section \ref{subsection: The effect on distance}.

\subsection{The balanced metrics}\label{subsection: The balanced metrics}

The metrics obtained by taking the coefficients of the polynomial $\alpha(1+cX)^k,$ i.e. $a_q = \alpha c^q{k \choose q},$ for any $\alpha, c>0,$ are fixed for both the $T$ and the $T_K$ maps; it is not for $T_{\nu}$ unless $c=1,$ in which we get the {\it round metric} -- the only palindromic balanced metrics for any $k$.  This can be explained by the fact that both the $T$ and $T_K$ maps respect the induced action of $\mathrm{SL}(2, \mathbb{C})$ on the space of metrics, while $T_{\nu}$ does not.

Starting with arbitrary $G = (a_0, a_1, \ldots, a_k)$ it is not entirely clear which balanced metric iterations of any of the operators $T, T_{\nu}, T_K$ will tend towards;  all we can say is the coefficients will be of the form $B = (b_0, \ldots, b_k)$ where $b_q = \alpha c^q {k \choose q}$ for some $\alpha, c >0$, and if $G$ is palindromic or the operator is $T_{\nu}$ then $c=1.$  We also note that when $k=2$ we can calculate the value $c$ as $c=\sqrt{a_2/a_0}$ and thus the balanced metric will be of the form $\alpha ( a_0, 2\sqrt{a_0a_2} , a_2 )$ for some scalar $\alpha>0$.

\subsection{Asymptotic behavior}\label{subsection: Asymptotic behavior}

In the long run the behavior of the iterations of $F = T, T_{\nu}, T_K$ is predictable.  For each function the limiting ratio
\[
    \sigma_{F, k} := \lim_{r \rightarrow \infty}
        \frac{
            \mathrm{dist}(F^{(r+1)}(G),\ F^{(\infty)}(G))
        }{
            \mathrm{dist}(F^{(r)}(G),\ F^{(\infty)}(G))
        }
\]
exists, and converges to a simple limit.  In \cite{donaldson5} Donaldson proves that in the case of the $T_{\nu}$ iteration and starting with a palindromic metric this $\sigma$-value can be computed as
\begin{equation}
    \sigma_{T_{\nu},k} =
        \frac{(k-1)k}{(k+2)(k+3)} \quad  \mathrm{(if\ }G\mathrm{\ is\ palindromic)}
\end{equation}
By examining many examples we also observed that if $G$ is not palindromic we get
\begin{equation}
    \sigma_{T_{\nu},k} =
        \frac{k}{k+2} \quad  \mathrm{(if\ }G\mathrm{\ is\ not\ palindromic)} \\
\end{equation}
while in the case of the $T$ iteration we have
\begin{equation}\label{sigmaT asymptotic}
    \sigma_{T,k} = \frac{(k-1)(k+6)}{(k+2)(k+3)}
\end{equation}
and for $T_K$ we get
\begin{equation}\label{sigmaTK asymptotic}
    \sigma_{T_K,k} = \frac{k-1}{k+3} .
\end{equation}
In neither of these latter two cases does it matter if we start with a palindromic metric or not.

We see that when $k=2$ we have
\[
    \sigma_{T_{\nu},2} \mathrm{(not\ pal.)} > \sigma_{T,2} > \sigma_{T_K,2} > \sigma_{T_{\nu},2} \mathrm{(pal.)}
\]
while for $k \geq 3$ we have
\[
    \sigma_{T,k} \geq \sigma_{T_{\nu},k} \mathrm{(not\ pal.)} > \sigma_{T_K,2} > \sigma_{T_{\nu},2} \mathrm{(pal.)}
\]
with strict inequalities for every $k>3.$  So in general, if we start with a palindromic metric $G$ we expect that the $T_{\nu}$ iterations will converge the most quickly, followed by $T_K$ and then by $T.$  Starting with a non-palindromic $G$ the $T_{\nu}$ iterations will slow down, and we find that $T_K$ will converge fastest.  Here $T$ is still slowest to converge.

\subsection{The effect on distance}\label{subsection: The effect on distance}

Despite this simple long-term behavior of the $T, T_{\nu},$ and $T_K$ iterations, the early behavior is still somewhat mysterious.  Perhaps one surprising fact along these lines is that {\it in general the $T$ operator is not distance reducing on the space of Hermitian metrics on $\mathrm{H}^0\left(\mathbb{CP}^1, O(k)\right)$}.  An example when $k=6$ is given in section \ref{section: Examples}.  This is the smallest value of $k$ for which the author has found such an example.

In \cite{cc} Calabi and Chen show that the Calabi flow is, in a certain sense, distance reducing.  Hence it might be surprising that $T$ is not given the expectation that it can be viewed as a discrete version of such a flow.

While it can happen that $T(G)$ is farther from the balanced metric than $G$ is, it does not appear to be the case that it can be {\it arbitrarily} farther.  Indeed for each of the operators $T, T_{\nu}, T_K$ the amount it can ``magnify'' the distance from the balanced metric appears to be simply bounded by a slow function of $k$.  This leads us to conjecture a bound for how far the $r$th iteration of any of the operators can be from the balanced metric.

Let $F = T, T_{\nu}, T_K$, let $G$ be any metric, and set $B = F^{(\infty)}(G)$ to be the balanced metric which the dynamical system $\{ F^{(r)}(G), r = 0, 1, 2, \ldots \}$ converges to.  Recall that we define
\[
    \mathrm{err}_{F, k} (G, r) = \mathrm{dist}\left( F^{(r)}(G),\ B \right).
\]
Let $d$ denote the initial distance from $G$ to $B$ in the space of Hermitian metrics.  Then we propose that in fact
\begin{equation}\label{bound}
    \mathrm{err}_{F, k} (G, r) < \log\left(1 + e^{kd} \sigma_{F, k}^r\right)
\end{equation}
for every $k>1$.  We do not expect this bound to be sharp.


\section{Examples}\label{section: Examples}

In this section we illustrate the findings from section \ref{section: Findings} with some examples.  We will always scale all metrics uniformly so that the limiting balanced metric begins with a one.  Note that each of the operators $T, T_{\nu}, T_K$ respect scaling.

We begin with $k=2$ and a non-palindromic metric proportional to $G = (1, 17, 36),$ and consider the $T_K$ iterations.  According to section \ref{subsection: The balanced metrics} the limiting balanced metric will be, after scaling, $B=(1, 12, 36).$  Below is a table displaying the results:  the first column gives the iteration $r$; the next three the entries of the metric; the second to last gives the distance from the balanced metric, or $\mathrm{err}_{T_K, k}(r, G)$; and the last column gives the bound $\mathrm{bnd}_{T_K,k}(d,r) = \log(1+e^{kd}\sigma_{T_K,k}^r)$ proposed in secion \ref{subsection: The effect on distance}.

\[
\begin{array}{c|ccc|c|c}
r & a_0 & a_1 & a_2 & \mathrm{dist}(-,B) & \mathrm{bnd} \\
\hline
  0  &     0.8826  &    15.0043  &    31.7738  &  0.2848  &  1.0180  \\
  1  &     0.9738  &    12.6377  &    35.0561  &  0.0640  &  0.3027  \\
  2  &     0.9946  &    12.1292  &    35.8067  &  0.0131  &  0.0683  \\
  3  &     0.9989  &    12.0259  &    35.9612  &  0.0026  &  0.0140  \\
  4  &     0.9998  &    12.0052  &    35.9922  &  0.0005  &  0.0028  \\
  5  &     1.0000  &    12.0010  &    35.9984  &  0.0001  &  0.0006
\end{array}
\]

We consider another non-palindromic metric, proportional to $G=(1, 25, 0.07, 13),$ with $k=3.$  We use the $T_{\nu}$ operator and list the results of the first few iterations $T_{\nu}^{(r)}(G)$ below.  We note that the limiting metric is $B = (1, 3, 3, 1).$

\[
\begin{array}{c|cccc|c|c}
r & a_0 & a_1 & a_2 & a_3 & \mathrm{dist}(-,B) & \mathrm{bnd} \\
\hline
  0  &    0.20720  &    5.18011  &    0.01450  &    2.69366  & 5.67338  & 17.02014  \\
  1  &    0.57206  &    2.68260  &    3.45522  &    1.58209  & 0.74488  & 16.50932  \\
  2  &    0.73295  &    2.72858  &    3.31411  &    1.32528  & 0.44129  & 15.99849  \\
  3  &    0.83372  &    2.82894  &    3.18320  &    1.18836  & 0.26423  & 15.48766  \\
  4  &    0.89777  &    2.89557  &    3.10812  &    1.11040  & 0.15845  & 14.97684  \\
  5  &    0.93773  &    2.93684  &    3.06435  &    1.06526  & 0.09505  & 14.46601  \\
 10  &    0.99505  &    2.99505  &    3.00496  &    1.00497  & 0.00739  & 11.91189  \\
 15  &    0.99961  &    2.99962  &    3.00039  &    1.00039  & 0.00057  & 9.35784  \\
 20  &    0.99997  &    2.99997  &    3.00003  &    1.00003  & 0.00004  & 6.80474
\end{array}
\]


We give one more table, this time beginning with a metric which moves away from the limiting metric after the first application of the operator $T.$  We choose the palindromic
\[
    G= (1,\ 6000,\ 150000,\ 20000000000,\ 150000,\ 6000,\ 1),
\]
with $k=6.$  Each iterate will be of the form
$(a_0, a_1, a_2, a_3, a_2, a_1, a_0),$ so we only keep track of $a_0,a_1, a_2, a_3.$  Again we uniformly scale so the limiting metric is exactly $B = (1, 6, 15, 20, 15, 6, 1).$
\[
\begin{array}{c|cccc|c|c}
r &  a_0 & a_1   &  a_2     &  a_3     &    \mathrm{err}     &     \mathrm{bnd}  \\
\hline\vspace{-10pt} \\
  0  &       0.00010  &       0.58903  &      14.72580  & 1963439.38600  & 17.69856  & 106.19139  \\
  1  &       0.00010  &       0.48814  &    1073.02459  &  733382.16850  & 18.10011  & 106.00906  \\
  2  &       0.00011  &       0.60722  &    1196.93120  &  414634.58830  & 17.67812  & 105.82674  \\
  3  &       0.00013  &       0.72695  &    1195.91914  &  257759.72070  & 17.21170  & 105.64441  \\
  4  &       0.00016  &       0.84269  &    1147.31003  &  167930.51810  & 16.72422  & 105.46208  \\
  5  &       0.00020  &       0.95726  &    1076.08572  &  112611.11230  & 16.22342  & 105.27976  \\
 10  &       0.00068  &       1.58083  &     669.18359  &   18910.93755  & 13.62571  & 104.36813  \\
 20  &       0.01002  &       3.32601  &     190.00391  &     970.58975  &  8.42894  & 102.54488  \\
 30  &       0.11205  &       5.07732  &      52.17933  &     117.34474  &  3.98456  & 100.72162  \\
 40  &       0.51092  &       5.88292  &      22.24884  &      34.20518  &  1.22538  & 98.89836  \\
 50  &       0.87358  &       5.99470  &      16.26035  &      22.28184  &  0.24744  & 97.07511  \\
 60  &       0.97741  &       5.99984  &      15.20684  &      20.36883  &  0.04187  & 95.25185  \\
 70  &       0.99629  &       6.00000  &      15.03350  &      20.05958  &  0.00682  & 93.42860  \\
 80  &       0.99940  &       6.00000  &      15.00541  &      20.00962  &  0.00110  & 91.60534  \\
 90  &       0.99990  &       6.00000  &      15.00088  &      20.00156  &  0.00018  & 89.78209  \\
100  &       0.99998  &       6.00000  &      15.00014  &      20.00025  &  0.00003  & 87.95883
\end{array}
\]

\begin{figure}[!p]
  \begin{center}
  \newcommand{\www}{100pt}
  \includegraphics[viewport=90 90 310 310, clip, width=1.5in]{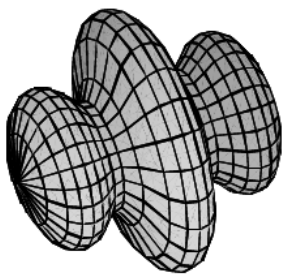}
  \end{center}
  \caption{$\mathbb{CP}^1$ with metric induced from $G = (1,300,300,300,1).$}
  \label{plot0}
\end{figure}

We finish the Riemann sphere case with a visual example of Donaldson's $T$-iterations.  We choose a palindromic metric which we can realize as induced from an embedding of $\mathbb{CP}^1$ into $\mathbb{R}^3.$  In particular we pick
\[
    G = (1, 300, 300, 300, 1)
\]
on $\mathrm{H}^0(\mathbb{CP}^1, O(4)),$ which is a metric obtained if one where to pinch the sphere around two latitudes giving it two narrow necks.  See figure \ref{plot0}.

Now, in figure \ref{The iterations} we plot the evolution of the metric $G$ under the iterations of $T, T_K$, and $T_{\nu},$ respectively.

\newcommand{\wid}{0.8in}
\begin{figure}[!p]
  \begin{center}
  \makebox[.3in]{}
  \makebox[\wid]{$r=0$}
  \makebox[\wid]{$r=1$}
  \makebox[\wid]{$r=2$}
  \makebox[\wid]{$r=3$}
  \makebox[\wid]{$r=4$}\\
  \raisebox{.3in}{\makebox[.3in]{$T$:}}
  \includegraphics[width=\wid]{plot0_gray.eps}
  \includegraphics[width=\wid]{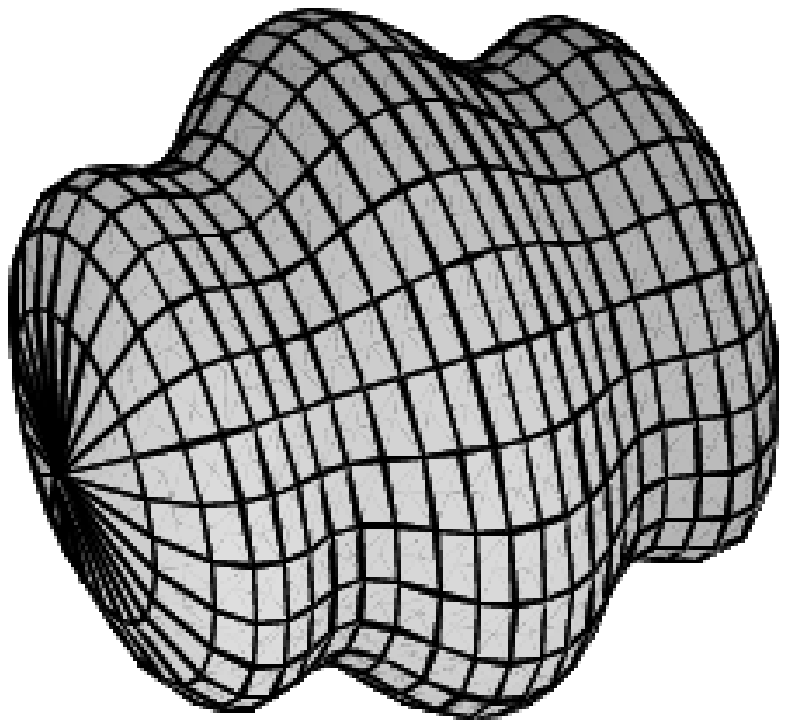}
  \includegraphics[width=\wid]{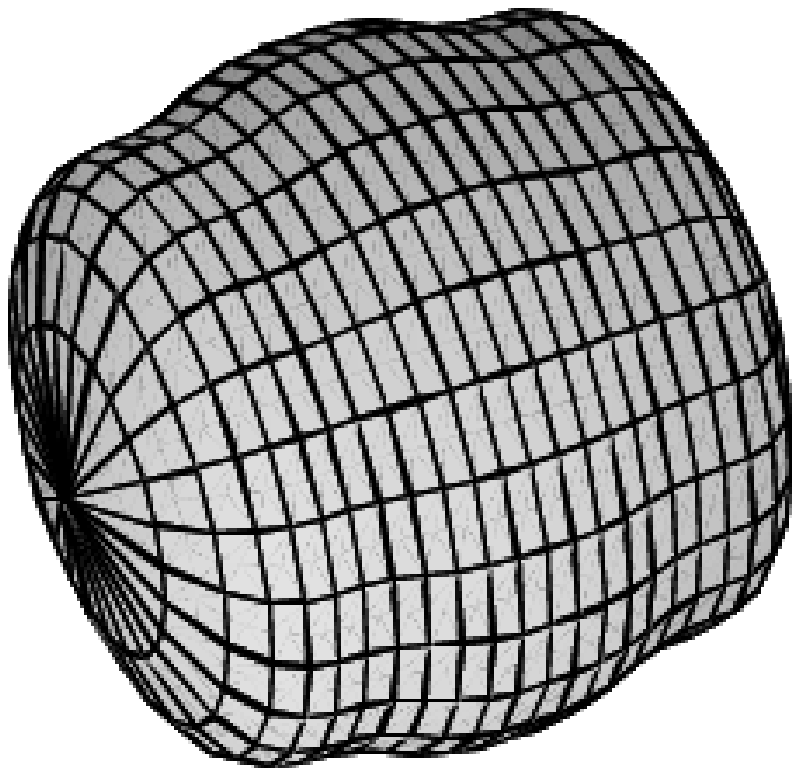}
  \includegraphics[width=\wid]{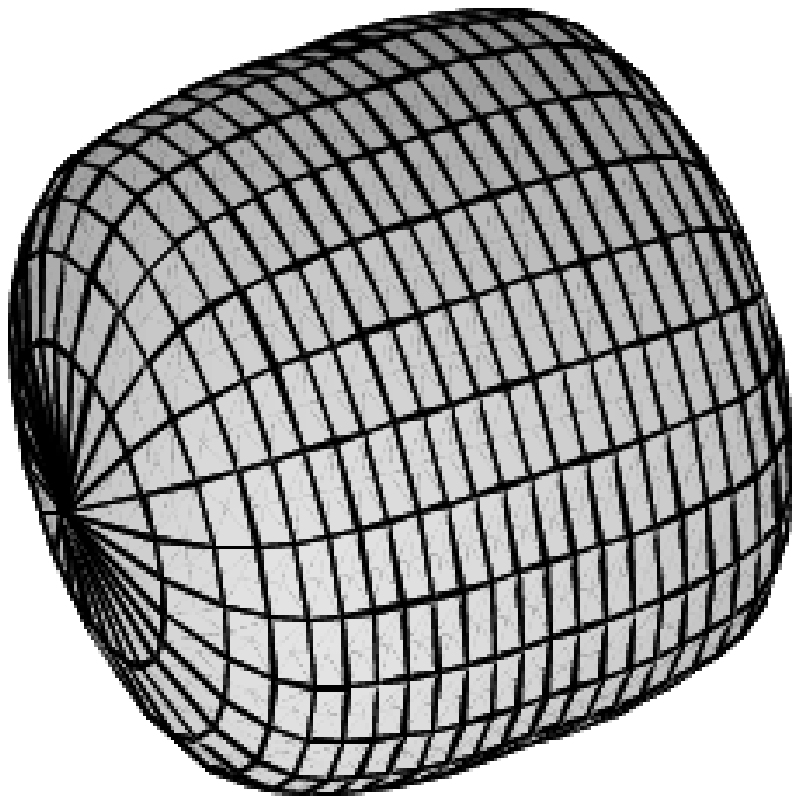}
  \includegraphics[width=\wid]{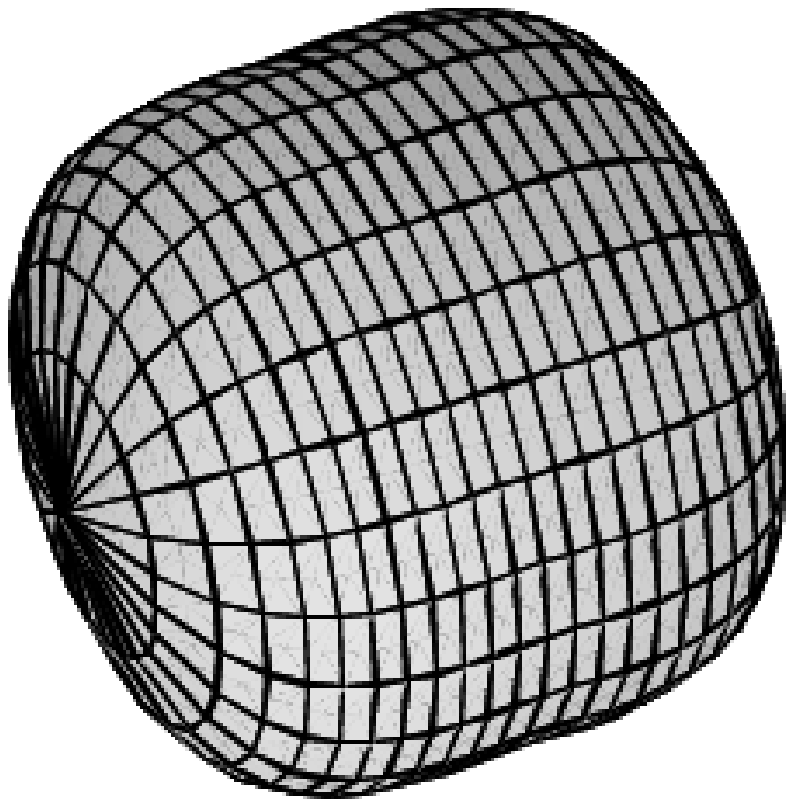}\\
  \raisebox{.3in}{\makebox[.3in]{$T_K$:}}
  \includegraphics[width=\wid]{plot0_gray.eps}
  \includegraphics[width=\wid]{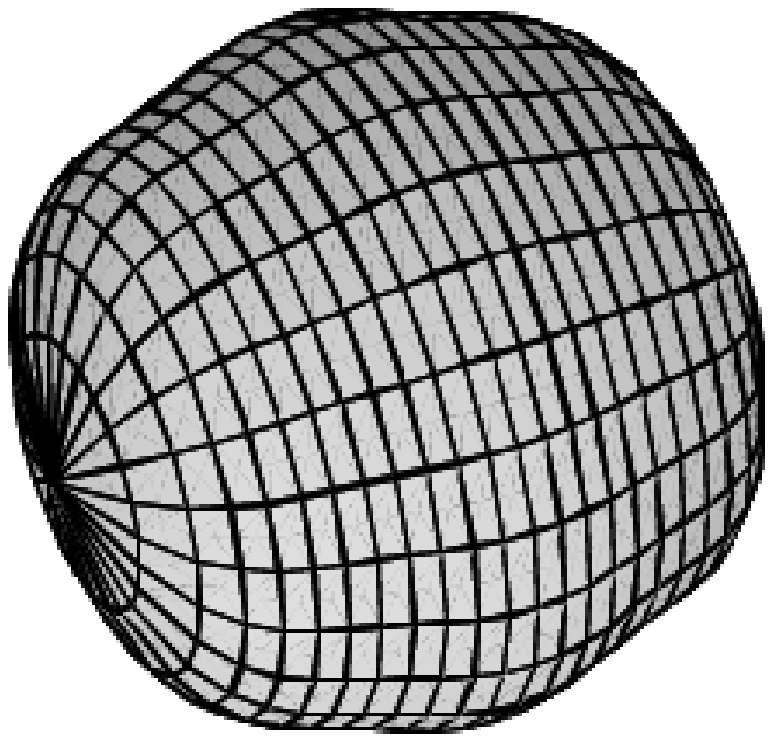}
  \includegraphics[width=\wid]{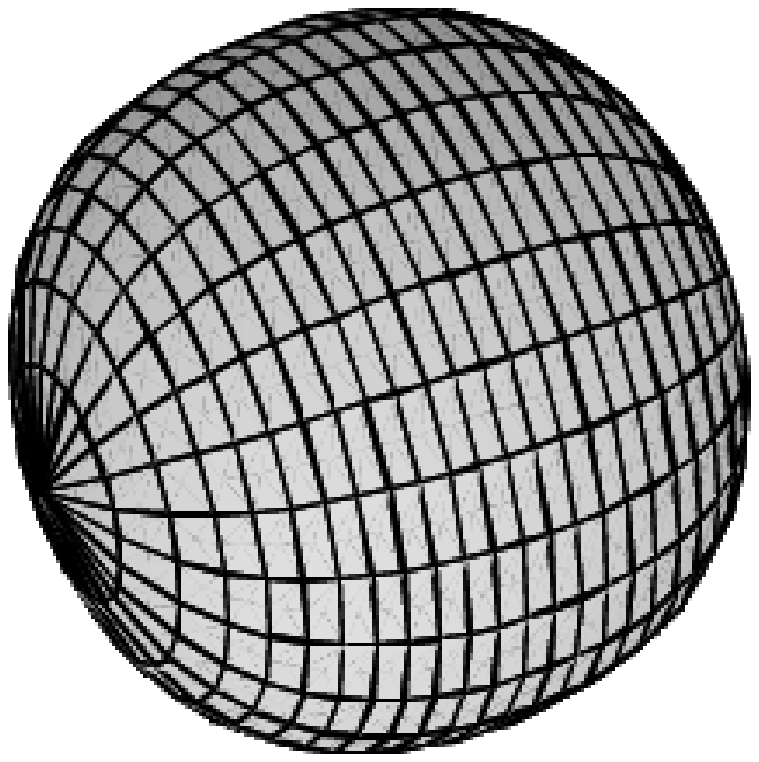}
  \includegraphics[width=\wid]{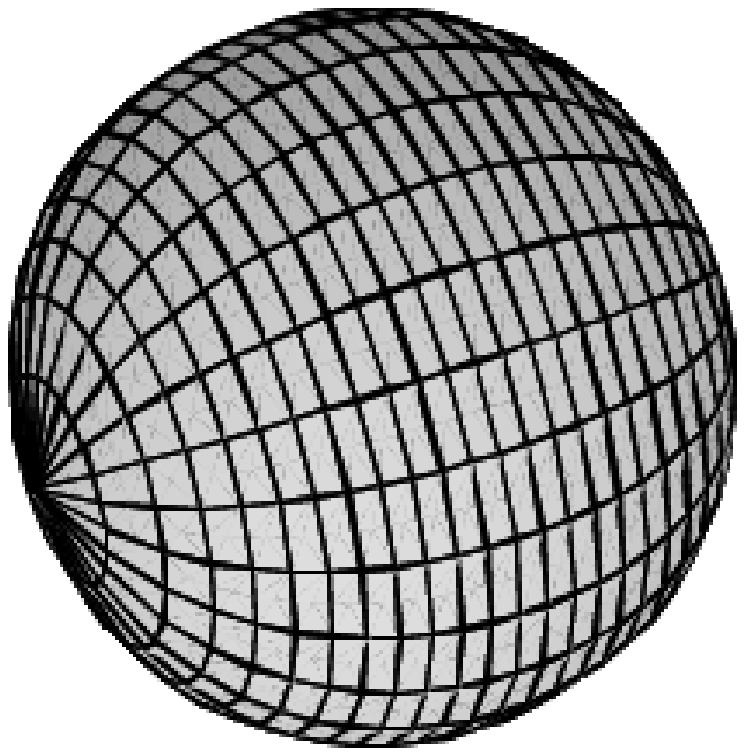}
  \includegraphics[width=\wid]{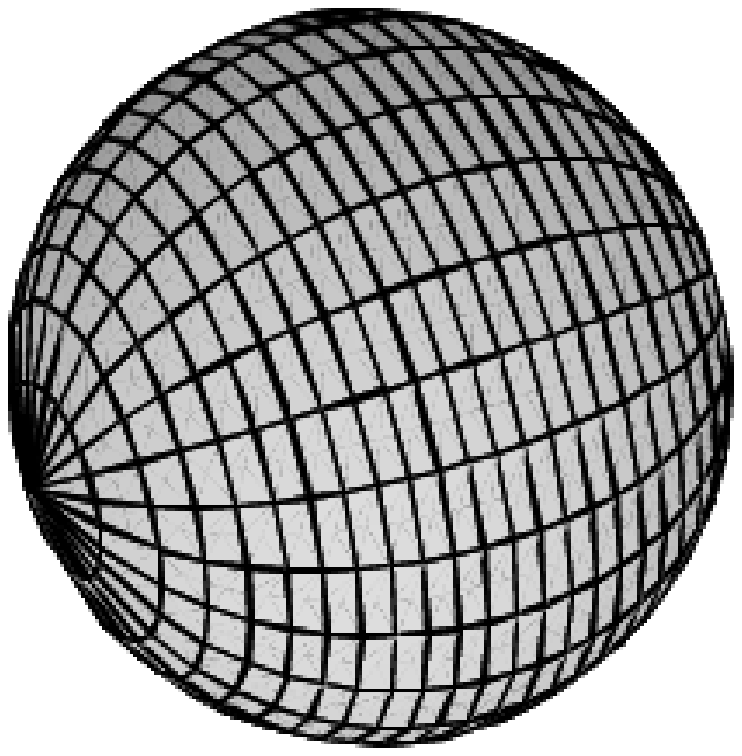}\\
  \raisebox{.3in}{\makebox[.3in]{$T_{\nu}$:}}
  \includegraphics[width=\wid]{plot0_gray.eps}
  \includegraphics[width=\wid]{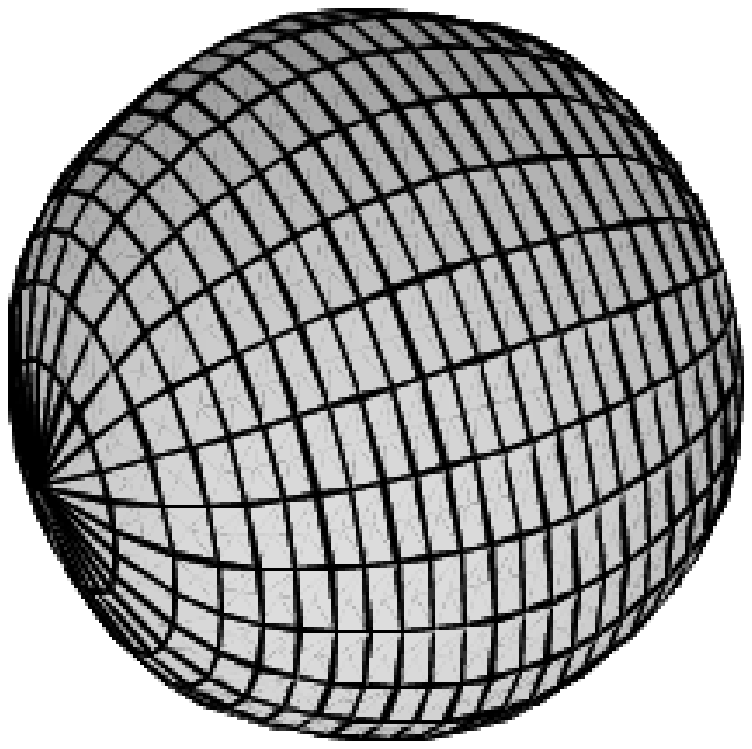}
  \includegraphics[width=\wid]{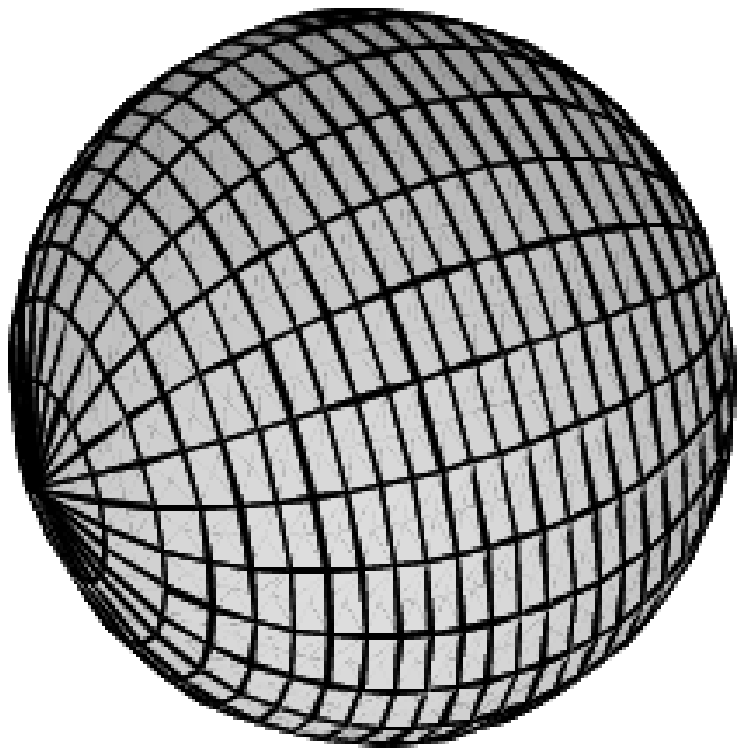}
  \includegraphics[width=\wid]{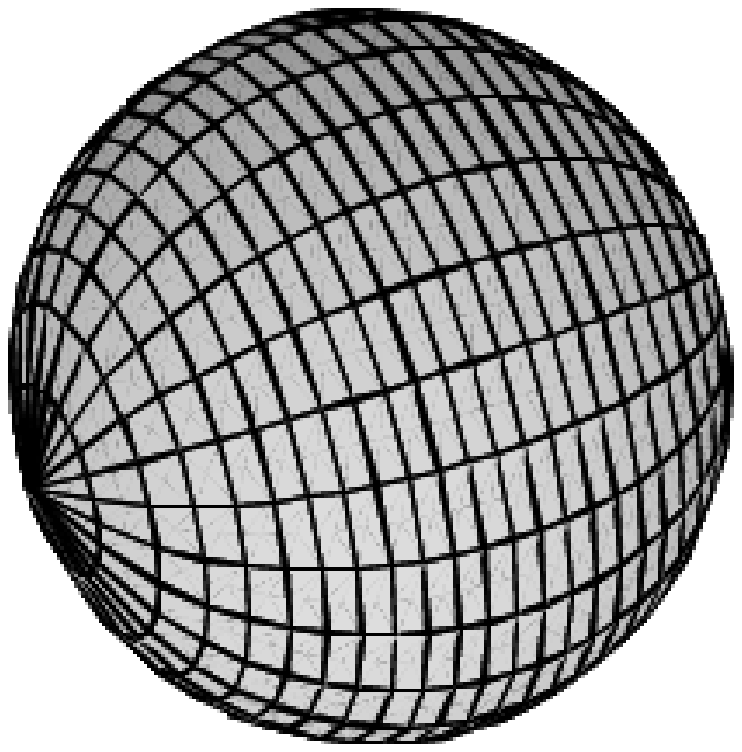}
  \includegraphics[width=\wid]{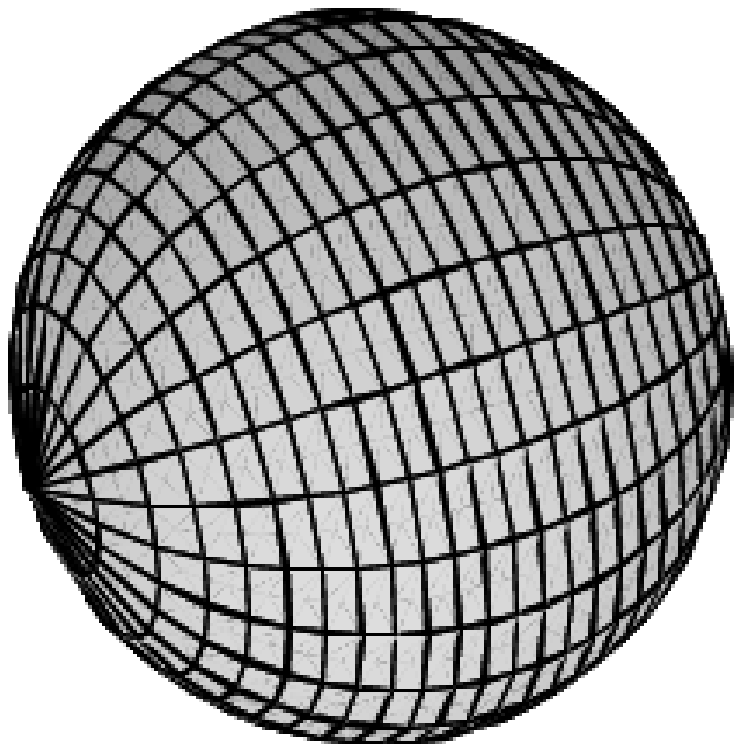}
  \end{center}
  \caption{The first four iterations.}
  \label{The iterations}
\end{figure}

Clearly the $T$ iterations are much slower in converging to a round sphere.  Not until the 3rd iteration does it become convex.  At the other extreme lie the $T_{\nu}$ iterations where the first iteration is already almost indistinguishable from a round sphere.  Intermediate between the two are the $T_K$ iterations.  This visually depicts the the observations in section \ref{subsection: Asymptotic behavior}, where rates of convergence were compared using asymptotic behavior.


\newcommand{\C}{\mathbb{C}}
\newcommand{\Proj}{\mathbb{P}}
\newcommand{\CP}{\mathbb{CP}}
\newcommand{\Tn}{T_{\nu}}

\section{Higher dimensional projective space}\label{section: Higher dimensional projective space}

Let us now investigate the complex projective space $X = \CP^n$ where $n>1$.  We will consider exclusively the $\Tn$ iteration.  Let $z_1, \ldots, z_n$ be local coordinates on $X=\CP^n$.  Let us fix once and for all a volume form $\nu$ on $X$ by using that induced by the normalized Fubini-Study metric. That is if
\begin{gather*}
    \omega = \frac{\sqrt{-1}}{2\pi}\ \partial \overline{\partial}\, \log \left( 1+ \sum | z_k |^2 \right)
    \phantom{XXXXX} \\[3mm]
        \phantom{XXXXX} =
        \frac{\sqrt{-1}}{2\pi}\
        \frac{
                \sum_{i,j} \left[ \left( 1 + \sum_k |z_k|^2 \right) \delta_{ij} - z_j \overline{z}_i \right]
                }{
                \left( 1 + \sum_k |z_k|^2 \right)^2
        }
        dz_i \wedge d\overline{z}_j
\end{gather*}
is the normalized Fubini-Study metric in local coordinates then we set
\[
    \nu = \omega^n =\, n!\, \left( \frac{\sqrt{-1}}{2\pi} \right)^n
        \frac{
            dz_1 \wedge d\overline{z}_1 \wedge \cdots \wedge dz_n \wedge d\overline{z}_n
            }{
            \left( 1 + |z_1|^2 + \ldots + |z_n|^2 \right)^{n+1}
        }.
\]
It is not hard to check that with this choice of volume form we get
\[
    \mathrm{Vol}(\CP^n) = \int_{\CP^n}\!\!\! \nu = 1.
\]

Again we set $L = O(1)$ and fix a $k>0$.  Note that a basis of $H^0(X,L^k)$ is given by the set of monomials in the $z_i$ of total degree $\leq k$.  Denote these by $w_1, \ldots, w_N$ where $N={n+k \choose k}$.  In this set-up we are studying embeddings
\[
    X = \CP^n \hookrightarrow \Proj \left( H^0(\CP^n, O(k)) \right) \cong \CP^{N-1}.
\]

As above we take $h$ to be the metric on $L^k = O(k)$ defined by
\[
    h = \left( \sum_{i=1}^N |w_i|^2 \right)^{-1}.
\]
Now if $G$ is a (positive definite Hermitian) matrix on $H^0(X,L^k)$ then $\Tn(G)$ is the matrix giving rise to the norm
\[
    \| s \|_{\mathrm{Hilb}_{\nu}} = R_{\nu} \int_X |s|_h^2 \nu,
\]
where
\[
    R_{\nu} = \frac{\mathrm{dim}\,  H^0(X, L^k)}{\mathrm{Vol}(X, \nu)} = N.
\]

The matrix $G$ has rows and columns indexed by the terms $w_i, \ i=1, \ldots, N$.  Let us take a diagonal matrix with terms $a_{i}^{-1}$.  Such a matrix corresponds to an (algebraic) metric invariant under the torus $\Lambda_n = (S^1)^n$ action $z_l \mapsto e^{i\theta_l}z_l,\ l=1,\ldots,n$.

An orthonormal basis, according to $G$, is given by
\[
    \left\{\sqrt{a_{i}}\,w_i, i=1, \ldots, N \right\}.
\]
Then in terms of the $a_i$'s the matrix $\Tn(G)$ will have diagonal entries $\tilde{a}_i^{-1}$ equal to
\[
    \| w_i \|_{\mathrm{Hilb}_{\nu}} =
    N n! \left( \frac{\sqrt{-1}}{2\pi} \right)^n \int_{\C^n} \frac{
        |w_i|^2 dz_1\wedge d\overline{z}_1\wedge\cdots \wedge dz\wedge d\overline{z}_n
        }{
        \left( \sum_{p=1}^N a_{p} |w_p|^2 \right) \left( 1+\sum_{q=1}^n |z_q|^2 \right)^{n+1}
        }.
\]
Changing to polar coordinates $z_j = r_j \exp\left(\sqrt{-1}\ \theta_j\right)$ and substituting $x_j=r_j^2$ we get
\begin{equation}\label{eqn: CPn integral}
    T(G)_{ii}=\tilde{a}_{i}^{-1} =
    N n! \int_0^{\infty}\hspace{-10pt}\cdots\!\!\int_0^{\infty}
    \frac{w_i(x) dx_1\cdots dx_n}{
        \left(\sum_{p=1}^N a_p w_p(x) \right) \left(1 + \sum_{q=1}^n x_q \right)^{n+1}
    }
\end{equation}
where $w(x)$ denotes the monomial $w$ with the substitutions $x_k = z_k, k=1, \ldots, n$.

\subsection{Asymptotic behavior in higher dimensions}

Let us consider the asymptotic behavior of $\Tn$.  Recall (see section \ref{subsection: Asymptotic behavior}) that in the case of $n=1$, i.e. when $X=\CP^1$, we defined
\[
    \sigma_{T_{\nu},k} := \lim_{r \rightarrow \infty}
        \frac{
            \mathrm{dist}(\Tn^{(r+1)}(G),\ \Tn^{(\infty)}(G))
        }{
            \mathrm{dist}(\Tn^{(r)}(G),\ \Tn^{(\infty)}(G))
        }.
\]
This value depends on whether or not the initial metric $G$ is invariant under the inversion map $z \mapsto z^{-1}$, or equivalently in homogeneous coordinates $Z_0 \leftrightarrow Z_1$.  In \cite{donaldson5} Donaldson computes these values theoretically, and our investigations corroborate his result:
\begin{equation}\label{eqn: sigma for CP1}
    \sigma_{T_{\nu},k} =
        \left\{ \begin{array}{rl}
            \frac{(k-1)k}{(k+2)(k+3)} & \mathrm{if\ }G\mathrm{\ is\ inv.\ under\ }Z_0 \leftrightarrow Z_1 \\
            \frac{k}{(k+2)} & \mathrm{otherwise}
        \end{array} \right. .
\end{equation}
Our goal is to show evidence for a simple extension of this formula valid on $X=\CP^n,\ n\geq 1$.

When $n>1$ there are many possible ways to extend the notion of a ``palindromic'' metric (as we defined in section \ref{section: The T operators}):  for the Riemann sphere we have those metrics invariant under $Z_0 \leftrightarrow Z_1$ but in general there are many permutations of the homogeneous coordinates $Z_0, \ldots, Z_n$ and it is trivial to check that if $G$ is a metric invariant under such a symmetry then so is $\Tn(G)$.  We might then expect that there can be distinct values for $\sigma$ depending on various symmetries the metric $G$ could be invariant under.  Thus we may find a different value for each (conjugacy class of) subgroup of $\mathrm{Sym}(n+1)$ -- the symmetric group on $n+1$ characters -- corresponding to metrics $G$ invariant under the automorphisms $Z_i \mapsto Z_{\pi(i)},\ i=0, 1, \ldots, n$ for $\pi$ ranging in the subgroup.

We present here some numerical findings in the cases of $n=2$ and $n=3$.  The itereated integrals (equation \ref{eqn: CPn integral}) grow in computational complexity quickly with increasing $n$.

We start with a metric $G$ which is torus-invariant, but otherwise `random' in the sense that it is not invariant under any permutation of the homogeneous coordinates.  We tabulate approximate numerical values for the asymptotic constant $\sigma$ here, all computed starting with `random' (but torus-invariant) metrics:

\begin{center}
\begin{tabular}{|c|c|c|c|c|}
    \hline
    $\sigma$    & $k=2$ & 3     & 4     & 5 \\
    \hline
    $n=2$       & 0.40  & 0.50  & 0.57    & 0.63 \\
    \hline
    3           & 0.33  & 0.43  & 0.50   &  0.56\\
    \hline
\end{tabular}
\end{center}

For the moment let us just note that the above values apparently follow the pattern:
\begin{equation}\label{eqn: sigma for Tv generic case}
    \sigma = \frac{k}{k+n+1}.
\end{equation}
When $n=1$, the fundamental case which we considered, this formula specializes to (\ref{eqn: sigma for CP1}).

The non-generic case, when $G$ might be invariant under a permutation of the homogeneous coordinate variables we find  simple behavior:
\begin{itemize}
    \item If there is no fixed-point-free permutation of the homogeneous coordinate variables under which $G$ is invariant then $\sigma$ is the same as computed in the asymmetric case.
    \item Otherwise suppose $G$ is invariant under some fixed-point-free permutation of the homogeneous coordinate variables.  Then we get new values for $\sigma$, as tabulated in the following table:
\end{itemize}

\begin{center}
\begin{tabular}{|c|c|c|c|c|}
    \hline
    $\sigma$    & $k=2$ & 3     & 4     & 5 \\
    \hline
    $n=2$       & 0.07 & 0.14 & 0.21 & 0.28 \\
    \hline
    3           & 0.05 & 0.11 & 0.17 & 0.22 \\
    \hline
\end{tabular}
\end{center}

One can check that approximate fractional equivalents to these numbers follow the pattern
\begin{equation}\label{eqn: sigma for general n in sym case}
    \sigma = \frac{(k-1)k}{(k+n+2)(k+n+3)}.
\end{equation}
We should stress that when $n=1$ equation (\ref{eqn: sigma for general n in sym case}), together with (\ref{eqn: sigma for Tv generic case}), specializes to (\ref{eqn: sigma for CP1}).  This together with various experimental evidence leads the author to ask the following:

\begin{question} Let $G$ be a torus-invariant metric arising from a matrix on $H^0(\CP^n, O_{\CP^n}(k))$, and let $B = \Tn^{\infty}(G)$ be the limiting balanced metric under the $\Tn$ iteration.  Define
\[
    \sigma_G(n,k) :=
        \lim_{r \rightarrow \infty}
            \frac{
                \mathrm{dist}(\Tn^{(r+1)}(G),\ B)
            }{
                \mathrm{dist}(\Tn^{(r)}(G),\ B)
            }.
\]
Let us say that $G$ is \emph{generally symmetric} if it is invariant under some fixed-point-free permutation of the homogeneous coordinates.  Then do we have the general formula
\begin{equation}\label{eqn:  sigma for Tv for general n}
    \sigma_G(n,k) =
        \left\{ \begin{array}{cl}
            \frac{(k-1)k}{(k+n+1)(k+n+2)} &
                \mathrm{if\ }G\mathrm{\ is\ generally\ symmetric} \\
            \frac{k}{(k+2)} & \mathrm{otherwise}
        \end{array} \right. \quad ?
\end{equation}
\end{question}

\subsection{Example computation}\label{subsection:  Example computation}

To illustrate a typical computation leading to some of the numbers above, take $n=3, k=4$.  Then
\[
    N = \dim  H^0\left(\CP^3, O(4)\right) = {3+4 \choose 4} = 35
\]
and a basis of $H^0(\CP^3, O(4))$ is (in local coordinates)
\begin{equation}\label{eqn: basis elements}
    \{w_1, \ldots, w_{35}\} = \{ 1, z_1, z_2, z_3, z_1^2, \ldots, z_2z_3^3, z_3^4\}.
\end{equation}
We choose a $G$ which is invariant under \emph{every} permutation of the homogeneous coordinates $Z_i,\, i=0, \ldots, 3$ (where $z_i=\frac{Z_i}{Z_0}$).  Taking into account these symmetries there are only five distinct basis elements:
\[
    1,\ z_1,\ z_1^2,\ z_1z_2,\ z_1z_2z_3.
\]
In the order the basis elements are listed in (\ref{eqn: basis elements}) -- first by degree then lexicographically -- these are the 1st, 2nd, 5th, 6th, and 15th elements.  In the notation used at the beginning of this section we pick diagonal entries of $G$: $a_i^{-1}$ in the row and column determined by the basis element $w_i$.  Due to the symmetries we will have five parameters:
\[
    G: a_1, a_2, a_5, a_6, a_{15}
\]
The iterations of $\Tn$ on these parameters -- denote them by $a_{i,r}, r=0, 1, \ldots, \infty$ -- will (after uniform scaling) tend toward the values $1, 4, 6, 12, 24$ for $i=1,2,5,6,15$ respectively.  This can readily be checked by noting the Fubini-Study metric is the balanced metric $B$.  At this point we should recall that the $a_i$ coefficients are actually entries in the inverse matrix $G^{-1}$, hence the entries of $G$ will tend to $1,1/4, 1/6, 1/12, 1/24$.  However it we can compute the approximate $\sigma$ values via
\[
    \lim_{r \rightarrow \infty}
            \frac{
                \mathrm{dist}(\Tn^{(r+1)}(G),\ B)
            }{
                \mathrm{dist}(\Tn^{(r)}(G),\ B)
            }.
    = \lim_{r \rightarrow \infty}
        \frac{
             \frac1{a_{2,r+1}} - \frac14
        }{
             \frac1{a_{2,r}} - \frac14
        }
    = \lim_{r \rightarrow \infty}
        \frac{a_{2,r+1}-4}{a_{2,r}-4}
\]
say (note the last equality follows since the $a_i,r$ are convergent).  Denote this last quotient, within the limit, as $\tilde{\sigma}_{r+1}$.  Its value should tend to the $\sigma$ value determining the asymptotic behavior of the $\Tn$ iterations on this metric.

Let us take $(a_1,a_2, a_5, a_6, a_{15}) = (1, 20, 30, 40, 50)$.  The limiting balanced metric will have corresponding coordinates proportional to $(1,4,6,12,24)$ as noted above.  However instead of uniformly scaling all metrics so the result is exactly this metric we will this time scale each metric so that its first coordinate (the $a_1$) is equal to one.  There is no loss of information: relation (\ref{eqn:  metric relation}) at the end of section (\ref{section: The T operators}) has the obvious adaptation to this situation; namely $\sum_{i=1}^{35} a_i/\tilde{a}_i = 35$.  Using this one can iteratively obtain the original numbers.  The advantage of doing this is that we no longer need to keep track of the first coordinate $a_1$.

With this convention we get the following table for the first eight iterations, as well as the approximate $\sigma$ values:
\[
\begin{array}{c|cccc|c}
\Tn^{(r)}& a_2 & a_5 & a_6 & a_{15} & \tilde{\sigma}_r \\
\hline
  0  &    20.0000000  &    30.0000000  &    40.0000000  &    50.0000000  &   0.0000  \\
  1  &     4.3071170  &     6.5967335  &    13.0915039  &    25.9850356  &   0.0192  \\
  2  &     4.0344368  &     6.0688663  &    12.1588436  &    24.3600437  &   0.1121  \\
  3  &     4.0052604  &     6.0105224  &    12.0258597  &    24.0613530  &   0.1528  \\
  4  &     4.0008611  &     6.0017223  &    12.0042908  &    24.0102741  &   0.1637  \\
  5  &     4.0001430  &     6.0002860  &    12.0007145  &    24.0017140  &   0.1661  \\
  6  &     4.0000238  &     6.0000476  &    12.0001191  &    24.0002857  &   0.1665  \\
  7  &     4.0000040  &     6.0000079  &    12.0000198  &    24.0000476  &   0.1666  \\
  8  &     4.0000007  &     6.0000013  &    12.0000033  &    24.0000079  &   0.1667

\end{array}
\]
We note that the apparent limiting value, $0.1\overline{6} = 1/6$, matches the value in equation (\ref{eqn:  sigma for Tv for general n}).


\section{Further Questions}

The case of a non-diagonal matrix (thus corresponding to a metric not invariant under $z\mapsto e^{i\theta}z$) was not treated in this paper.  Investigating this direction one might see whether the asymptotic values (see equations \ref{sigmaT asymptotic} and \ref{sigmaTK asymptotic} in section \ref{subsection: Asymptotic behavior}, or \ref{eqn:  sigma for Tv for general n} in section \ref{subsection:  Example computation}) remain valid, and whether the bound (\ref{bound}) given in section \ref{subsection: The effect on distance} still holds.  If the bound does still hold then it would be interesting to work towards a sharp bound.

In another direction one might ask whether or not the operators $T, T_{\nu}, T_K$ are distance decreasing {\it after} the first iteration; or put another way:  is the square of each of these operators distance reducing?  No counter example to this was found.

The next step is to look beyond $\CP^n$, perhaps to toric varieties (see for example \cite{bd}), $K3$ surfaces, Calabi-Yau $3$-folds, etc., and work out the same convergence properties of these dynamical systems.  It would also be interesting to compare the convergence properties of the $T$-iterations to those of PDE methods for finding canonical metrics, such as the Ricci flow.  All these questions the author hopes to examine later.

\bibliography{refs}
\bibliographystyle{alpha}

\end{document}